\documentclass[a4paper,11pt]{article}

\usepackage{a4wide}
\usepackage{amssymb}
\usepackage{amsmath}
\usepackage{color}
\usepackage{url}
\usepackage{hyperref}
\usepackage{algorithm,algorithmic}

\newcommand{\Rstar}[2]{\mathbb{R}_*^{#1\times #2}}
\newcommand{\St}[2]{\mathrm{St}(#1,#2)}

\def \D {{\mathrm{D}}}

\def \grad {{\mathrm{grad}}}

\def \trace {{\mathrm{trace}}}

\def \sym {{\mathrm{sym}}}
\def \skew {{\mathrm{skew}}}

\def \T {\mathrm{T}}

\def \GL {{\mathrm{GL}}}
\def \StX {{\mathrm{St}}}

\def \Exp {{\mathrm{Exp}}}

\def \Mp {{\mathcal{M}(p,m\times n)}}
\def \Rp {{\mathbb{R}_*^{m\times p}\times \mathbb{R}_*^{n\times p}}}
\def \TRp {{\mathbb{R}^{m\times p}\times \mathbb{R}^{n\times p}}}
\def \SRp {{\mathrm{St}(p,m)\times \mathbb{R}_*^{n\times p}}}

\def \dM {\dot{M}}
\def \dN {\dot{N}}
\def \dR {\dot{R}}
\def \dX {{\dot{X}}}
\def \dY {{\dot{Y}}}
\def \cM {\check{M}}
\def \cN {\check{N}}
\def \cX {\check{X}}
\def \rh {\mathrm{h}}
\def \rM {\mathrm{M}}
\def \rN {\mathrm{N}}

\newtheorem{thrm}{Theorem}[section]

\newtheorem{prpstn}[thrm]{Proposition}

\begin{document}

\title{Two Newton methods on the manifold of fixed-rank matrices endowed with Riemannian quotient geometries
  \footnotemark[1]} 
\author{P.-A. Absil\footnotemark[2] \and Luca Amodei\footnotemark[3] \and Gilles
  Meyer\footnotemark[4]}

\date{\today}

\maketitle

\renewcommand{\thefootnote}{\fnsymbol{footnote}} 
\footnotetext[1]{This
paper presents research results of the Belgian Network DYSCO
(Dynamical Systems, Control, and Optimization), funded by the
Interuniversity Attraction Poles Programme, initiated by the Belgian
State, Science Policy Office. The scientific responsibility rests with
its authors.}
\footnotetext[2]{Department of Mathematical Engineering, ICTEAM Institute, Universit\'e catholique de Louvain, B-1348 Louvain-la-Neuve, Belgium
(\url{http://sites.uclouvain.be/absil/}).}
\footnotetext[3]{Institut de Mathématiques de Toulouse, Université
  Paul Sabatier, 118 route de Narbonne, 31062 Toulouse cedex 9, France (\texttt{luca.amodei@math.univ-toulouse.fr}).}
\footnotetext[4]{Department of Electrical Engineering and Computer Science,
University of Li\`ege,
B-4000 Li\`ege, Belgium
(\texttt{gillesmy@gmail.com}).}
\renewcommand{\thefootnote}{\arabic{footnote}}

\begin{center}
{\bf Abstract}\\[1em]
\begin{minipage}{12cm}
\noindent We consider two Riemannian geometries for the manifold $\Mp$
of all $m\times n$ matrices of rank $p$. The geometries are induced on
$\Mp$ by viewing it as the base manifold of the submersion
$\pi:(M,N)\mapsto MN^\T$, selecting an adequate Riemannian metric on
the total space, and turning $\pi$ into a Riemannian submersion. The
theory of Riemannian submersions, an important tool in Riemannian
geometry, makes it possible to obtain expressions for fundamental
geometric objects on $\Mp$ and to formulate the Riemannian
Newton methods on $\Mp$ induced by these two geometries. The Riemannian
Newton methods admit a stronger and more streamlined convergence
analysis than the Euclidean counterpart, and the computational
overhead due to the Riemannian geometric machinery is shown to be
mild. Potential applications include low-rank matrix completion and
other low-rank matrix approximation problems. 
\end{minipage}
\end{center}


\section{Introduction}
\label{sec:intro}

Let $m$, $n$, and $p\leq\min\{m,n\}$ be positive integers and let
$\Mp$ denote the set of all rank-$p$ matrices of size $m\times n$,
\begin{equation}  \label{eq:Mp}
\Mp = \{X\in\mathbb{R}^{m\times n}: \mathrm{rank}(X)=p\}.
\end{equation}
Given a smooth function $f:\Mp\to\mathbb{R}$, we
consider the problem
\begin{equation}  \label{eq:prob-f}
\min f(X) \quad \text{subject to $X\in\Mp$}.
\end{equation}

Problem~\eqref{eq:prob-f} subsumes low-rank matrix approximation problems,
where $f(X)\equiv \|A-X\|^2$ with $A\in\mathbb{R}^{m\times n}$ given
and $\|\cdot\|$ a (semi)norm. In particular, it includes
low-rank matrix completion problems, which have been the topic of much
attention recently; see~\cite{KesMonOh2010,DaiMilKer2011,BouAbs2011,eyc2011,MisMeyBacSep2011,DaiKerMil2012} and
references therein. Interestingly, low-rank matrix
completion problems combine two sparsity aspects: only a few elements
of $A$ are available, and the vector of singular values of $X$ is
restricted to have only a few nonzero elements.

This paper belongs to a trend of research, see~\cite{HM94,HelSha1995,SimEld2010,eyc2011,MisMeySep2011,MisMeyBacSep2011}, where problem~\eqref{eq:prob-f} is tackled using differential-geometric techniques exploiting the fact that $\Mp$ is a submanifold of $\mathbb{R}^{m\times n}$. 
We are interested in Riemannian Newton methods (see~\cite{Smi94,ADM2002,AMS2008})
for problem~\eqref{eq:prob-f}, with a preference for the pure Riemannian setting~\cite{Smi94}. This setting involves defining a Riemannian metric on $\Mp$ and providing an expression for the Riemannian connection---which underlies the Riemannian Hessian---and for the Riemannian exponential. When $\Mp$ is viewed as a Riemannian submanifold of $\mathbb{R}^{m\times n}$, the necessary ingredients for computing the Riemannian Hessian are available~\cite[\S 2.3]{eyc2011}, but a closed-form expression of the Riemannian exponential has been elusive in that geometry.

In this paper, we follow a different approach that strongly relies on two-term factorizations of low-rank matrices. To this end, let
\begin{equation}  \label{eq:Rstar}
\Rstar{m}{p} = \{X\in\mathbb{R}^{m\times p}: \mathrm{rank}(X)=p\}
\end{equation}
denote the set of all full-rank $m\times p$ matrices, and observe that, since the function 
\begin{equation}  \label{eq:pi}
\pi: \Rp \to \Mp: (M,N)\mapsto MN^\T
\end{equation}
is surjective, problem~\eqref{eq:prob-f} amounts to
the optimization over its domain of the function $\bar{f}=f\circ\pi$, i.e.,
\begin{equation}  \label{eq:bar-f}
\bar{f}: \Rp\to\mathbb{R}: (M,N)\mapsto f(MN^\T).
\end{equation}
Pleasantly, whereas $\Mp$ is a nonlinear space, $\Rp$ is an open subset of a linear space; more precisely, $\Rp$ is the linear space $\TRp$ with a nowhere dense set excerpted. 
The downside is that the minimizers of $\bar{f}$ are never
isolated; indeed, for all $(M,N)\in\Rp$, $\bar{f}=f\circ\pi$ assumes the same value $\bar{f}(M,N)$ at all points of
\begin{equation}  \label{eq:pi-1}
\pi^{-1}(MN^\T) = \{(MR,NR^{-\T}): R\in\GL(p)\},
\end{equation}
where 
\[
\GL(p) = \{R\in\mathbb{R}^{p\times p}: \det(R)\neq0\}
\]
denotes the general linear group of degree $p$. 
In the context of Newton-type methods, this can be a source of concern since, whereas the convergence theory of
Newton's method to nondegenerate minimizers is well understood (see,
e.g.,~\cite[Theorem~5.2.1]{DS83}), the situation becomes
more intricate in the presence of non-isolated minimizers (see,
e.g.,~\cite{Gri1985}). 

The proposed remedy to this downside consists in elaborating a Riemannian Newton method that evolves conceptually on $\Mp$---avoiding the structural degeneracy in $\Rp$---while still being formulated in $\Rp$. This is made possible by endowing $\Rp$ and $\Mp$ with Riemannian metrics that turn $\pi$
into a Riemannian submersion. The theory of Riemannian submersions~\cite{ONe1966,ONe1983} then provides a way of representing the Riemannian connection and the Riemannian exponential of $\Mp$ in terms of the same objects of $\Rp$. 

It should be pointed out that the local quadratic convergence of the
Riemannian Newton method is retained if the Riemannian connection is
replaced by any affine connection and the Riemannian exponential is
replaced by any first-order approximation, termed retraction;
see~\cite[\S 6.3]{AMS2008}. The preference for the pure Riemannian
setting is thus mainly motivated by the mathematical elegance of a
method fully determined by the sole Riemannian metric.

Some of the material of this paper is inspired from the PhD thesis~\cite{Mey2011} and the talk~\cite{AmoDedYak2009}.

The paper is organized as follows. In the short
sections~\ref{sec:quot} and~\ref{sec:RS}, we show that $\pi$ is a
submersion and we recall some fundamentals of Riemannian submersions. A
first, natural but unsuccessful attempt at turning $\pi$ into a
Riemannian submersion is presented in Section~\ref{sec:nR-no}. Two
ways of achieving success are then presented in
sections~\ref{sec:R-no} and~\ref{sec:R-o}. In Section~\ref{sec:R-no},
the strategy consists of introducing a non-Euclidean Riemannian metric
on $\Rp$, whereas in Section~\ref{sec:R-o}, the plan of action is to restrict $\Rp$ by imposing orthonormality of one of the factors. We obtain closed-form expressions for the Riemannian connection (in both cases) and for the Riemannian exponential (in the latter case). Conclusions are drawn in Section~\ref{sec:conc}.

\section{$\Mp$ as a quotient manifold}
\label{sec:quot}

The set $\Mp$ of rank-$p$ matrices of size $m\times n$ is known to be an
embedded submanifold of dimension $p(m+n-p)$ of $\mathbb{R}^{m\times
  n}$, connected whenever $\max\{m,n\}>1$; see~\cite[Ch.~5, Prop.~1.14]{HM94}. Hence $\pi$~\eqref{eq:pi}
is a smooth surjective map between two manifolds. 

We show that $\pi$ is a submersion, i.e., that the differential of
$\pi$ is everywhere surjective. Observe that the tangent space to
$\Rp$ at $(M,N)$ is given by
\[
\T_{(M,N)}\Rp = \TRp;
\] 
this
comes from the fact that $\Rp$ is an open submanifold of the Euclidean
space $\TRp$~\cite[\S 3.5.1]{AMS2008}. For all $(M,N)\in\Rp$ and all
$(\dot{M},\dot{N})\in\TRp$, we have
$\mathrm{D}\pi(M,N)[(\dot{M},\dot{N})] = \dot{M}N^\T +
M\dot{N}^\T$. Working in a coordinate system where $M = \begin{bmatrix}
  I & 0 \end{bmatrix}^\T$ and $N = \begin{bmatrix}
  I & 0 \end{bmatrix}^\T$, one readily sees that the dimension of the
range of the map $(\dot{M},\dot{N}) \mapsto
\mathrm{D}\pi(M,N)[(\dot{M},\dot{N})]$ is equal to $p(m+n-p)$, the
dimension of the codomain of $\pi$. Hence $\pi$ is a submersion.

As a consequence, by the submersion
theorem~\cite[Proposition~3.3.3]{AMS2008}, the fibers $\pi^{-1}(MN^\T)$
are $p^2$-dimensional submanifolds of $\Rp$. Moreover,
by~\cite[Proposition~3.5.23]{AbrMarRat1988}, the equivalence relation
$\sim$ on $\Rp$, defined by $(M_a,N_a)\sim(M_b,N_b)$ if and only if
$\pi(M_a,N_a)=\pi(M_b,N_b)$, is regular and $\Rp/\sim$ is a quotient
manifold diffeomorphic to $\Mp$.

\section{Riemannian submersion: principles}
\label{sec:RS}

Turning $\pi$ into a Riemannian submersion amounts to endowing its
domain $\Rp$ with a Riemannian metric $\bar{g}$ that satisfies a certain
invariance condition, described next.

By definition, the vertical space $\mathcal{V}_{(M,N)}$ at a point $(M,N)\in\Rp$ is the tangent space to
the fiber $\pi^{-1}(MN^\T)$~\eqref{eq:pi-1}. We obtain
\begin{equation} \label{eq:V}
\mathcal{V}_{(M,N)} = \{(M\dot R,-N\dot R^\T): \dot
R\in\mathbb{R}^{p\times p}\}.
\end{equation}

Let $\bar{g}$ be a Riemannian metric on $\Rp$. Then one defines
the \emph{horizontal space} $\mathcal{H}_{(M,N)}$ at $(M,N)$ to
be the orthogonal complement of $\mathcal{V}_{(M,N)}$ in $\TRp$ relative to
$\bar{g}_{(M,N)}$, i.e.,
\begin{equation}  \label{eq:H-gen}
\mathcal{H}_{(M,N)} = \{(\dM,\dN)\in\TRp:
\bar{g}_{(M,N)}((\dM,\dN),(M\dR,-N\dR^\T))=0, \forall
\dR\in\mathbb{R}^{p\times p}\}.
\end{equation}
Next, given a tangent vector $\dot{X}_{MN^\T}\in \T_{MN^\T}\Mp$,
there is one and only one
\begin{equation}  \label{eq:lift-gen} 
\dot{X}_{(M,N)}\in\mathcal{H}_{(M,N)} 
\quad\text{such that}\quad \D\pi(M,N)[\dot{X}_{(M,N)}] = \dot{X}_{MN^\T},
\end{equation}
where $\D\pi(X)[\dX]$ denotes the differential of $\pi$ at $X$ applied to
$\dX$. 
This $\dot{X}_{(M,N)}$ is termed the \emph{horizontal lift} of
$\dot{X}_{MN^\T}$ at $(M,N)$.   (In order to lighten the notation, we use the
same symbol for a tangent vector to $\Mp$ and its horizontal lift; the distinction is clear from the subscript or from the context.) 
If (and only if), for all $(M,N)\in\Rp$, all
$\dX_{MN^\T},\cX_{MN^\T}\in \T_{MN^\T}\Mp$, and all
$R\in\GL(p)$, it holds that
\begin{equation}  \label{eq:g=g}
\bar{g}_{(M,N)}(\dX_{(M,N)},\cX_{(M,N)}) =
\bar{g}_{(MR,NR^{-\T})}(\dX_{(MR,NR^{-\T})},\cX_{(MR,NR^{-\T})}),
\end{equation}
then there is a (unique) Riemannian metric $g$ on $\Mp$ consistently defined by
\[
g_{MN^\T}(\dX_{MN^\T},\cX_{MN^\T}) =
\bar{g}_{(M,N)}(\dX_{(M,N)},\cX_{(M,N)}).
\]
The submersion $\pi:(\Rp,\bar{g})\to(\Mp,g)$ is then termed a \emph{Riemannian
  submersion}, and $(\Mp,g)$ is termed a \emph{Riemannian quotient
  manifold} of $(\Rp,\bar{g})$. (We will sometimes omit the
Riemannian metrics in the notation when they are clear from the context or undefined.)

In summary, in order to turn $\pi$ into a Riemannian submersion, we
``just'' have to choose a Riemannian metric $\bar{g}$ of $\Rp$
that satisfies the invariance condition~\eqref{eq:g=g}. 

\section{$\Mp$ as a non-Riemannian quotient manifold}
\label{sec:nR-no}

In this section, we consider on $\Rp$ the Euclidean metric $\bar{g}$,
defined by
\begin{equation}  \label{eq:Euc-bar-g}
\bar{g}_{(M,N)}\left( (\dM,\dN), (\cM,\cN) \right) := 
\trace(\dM^\T\cM) + \trace(\dN^\T\cN),
\end{equation}
and we show that the invariance
condition~\eqref{eq:g=g} does not hold. Hence $\pi:(\Rp,\bar{g})\to\Mp$
cannot be turned into a Riemannian submersion. 

The horizontal space~\eqref{eq:H-gen} is
\[
\mathcal{H}_{(M,N)} = \{(\dot{M},\dot{N}): \trace(\dM^\T M\dR)+\trace(-\dN^\T N\dR^\T) = 0, \forall
\dot{R}\in\mathbb{R}^{p\times p}\}.
\]
Using the identities $\trace(A)=\trace(A^\T)$ and
$\trace(AB)=\trace(BA)$, we obtain the identity
$\trace(\dM^\T M\dR)+\trace(-\dN^\T N\dR^\T) = 
\trace\left((\dR^\T (M^\T\dot{M} - \dot{N}^\T N)\right)$. It follows that the
following propositions are equivalent:
\begin{gather*}
(\dot{M},\dot{N})\in\mathcal{H}_{(M,N)},
\\  M^\T\dot{M} = \dot{N}^\T N,
\\ \exists L_M,L_N,S:
\begin{cases}
\dot{M} = M_\perp L_M + M(M^\T M)^{-1}S
\\ \dot{N} = N_\perp L_N + N(N^\T N)^{-1} S^\T,
\end{cases}
\end{gather*}
where $M_\perp$ denotes an orthonormal $m\times(m-p)$ matrix such that
$M^\T M_\perp=0$, and likewise for $N_\perp$. 

Let $X=MN^\T$ and let $\dX_{MN^\T}\in\T_{MN^\T}\Mp$. We seek an
expression for the \emph{horizontal lift} $\dX_{(M,N)}=(\dX_{\rM(M,N)},\dX_{\rN(M,N)})$ of $\dX_{MN^\T}$ at $(M,N)$, defined
by~\eqref{eq:lift-gen}.
By a reasoning similar to the one detailed in
Section~\ref{sec:R-no-hor-lift} below, we obtain
\[
\dX_{\rM(M,N)} = (\dX_{MN^\T} N-MK)(N^\T N)^{-1} \quad\text{and}\quad
\dX_{\rN(M,N)} = (\dX_{MN^\T}^\T M-NK^\T)(M^\T M)^{-1},
\]
where $K$ solves the Sylvester equation
\[
M^\T MK + KN^\T N = M^\T\dX_{MN^\T} N.
\]

One sees by inspection, or by a numerical check, that the invariance
condition~\eqref{eq:g=g} does not hold, and this concludes the
argument.

\section{$\Mp$ as a Riemannian quotient manifold of $\Rp$}
\label{sec:R-no}  

In this section, we proceed as in Section~\ref{sec:nR-no}, but now
with a different Riemannian metric $\bar{g}$, defined
in~\eqref{eq:bar-g} below. As we will see, the
rationale laid out in Section~\ref{sec:nR-no} now leads to the conclusion
that $\pi:(\Rp,\bar{g})\to\Mp$, with $\bar{g}$ given
by~\eqref{eq:bar-g} instead of~\eqref{eq:Euc-bar-g}, \emph{can} be turned into a Riemannian
submersion. This endows $\Mp$ with a Riemannian metric, $g$. We then work
out formulas for the Riemannian gradient and Hessian of $f$ on the Riemannian
manifold $(\Mp,g)$, and we state the corresponding Newton method.

\subsection{Riemannian metric in total space}

Inspired from the case of the Grassmann manifold viewed as a
Riemannian quotient manifold of $\mathbb{R}_*^{n\times
  p}$~\cite[Example~3.6.4]{AMS2008}, we consider the Riemannian metric
$\bar{g}$ on $\Rp$ defined by
\begin{equation}  \label{eq:bar-g}
\bar{g}_{(M,N)}\left( (\dM,\dN), (\cM,\cN) \right) := \trace
\left((M^\T M)^{-1} \dM^\T\cM + (N^\T N)^{-1} \dN^\T\cN \right).
\end{equation}
We now proceed to show that it satisfies the invariance
condition~\eqref{eq:g=g}. 

\subsection{Horizontal space}

The elements $(\dM,\dN)$ of the horizontal 
space $\mathcal{H}_{(M,N)}$~\eqref{eq:H-gen} are readily found to be
characterized by
\begin{equation}  \label{eq:H-cond}
M^\T\dM(M^\T M)^{-1} = (N^\T N)^{-1}\dN^\T N.
\end{equation}
In other words, 
\begin{equation}  \label{eq:H-def}
\mathcal{H}_{(M,N)} = \{(\dM,\dN)\in\TRp : N^\T NM^\T\dM = \dN^\T NM^\T M
\}.
\end{equation}

\subsection{Horizontal lift}
\label{sec:R-no-hor-lift}

Let $X=MN^\T$ and let $\dX_{MN^\T}$ belong to $\T_{MN^\T}\Mp$. 
We seek an
expression for the \emph{horizontal lift}
$\dX_{(M,N)}=(\dX_{\rM(M,N)},\dX_{\rN(M,N)})$ defined in~\eqref{eq:lift-gen}.
In view of~\eqref{eq:H-cond}, we find that the horizontality condition
$(\dX_{\rM(M,N)},\dX_{\rN(M,N)})\in\mathcal{H}_{(M,N)}$ is equivalent to 
\begin{subequations}  \label{eq:dXMN=LK}
\begin{align}
\dX_{\rM(M,N)} &= M_\perp L_\rM + M(M^\T M)^{-1}K(M^\T M)
\\ \dX_{\rN(M,N)} &= N_\perp L_\rN + N(N^\T N)^{-1}K^\T(N^\T N),
\end{align}
\end{subequations}
where 
$L_\rM\in\mathbb{R}^{(m-p)\times p}$,
$L_\rN\in\mathbb{R}^{(n-p)\times p}$ and $K\in\mathbb{R}^{p\times
  p}$. Since $\mathrm{D}\pi(M,N)[\dX_{\rM(M,N)},\dX_{\rN(M,N)}] \equiv M\dX_{\rN(M,N)}^\T +
\dX_{\rM(M,N)} N^\T$, the definition~\eqref{eq:lift-gen} implies that
\begin{equation}  \label{eq:dX}
\dX_{MN^\T} = M\dX_{\rN(M,N)}^\T + \dX_{\rM(M,N)} N^\T.
\end{equation}
Replacing~\eqref{eq:dXMN=LK} in~\eqref{eq:dX} yields
\begin{equation}  \label{eq:dX=LK}
\dX_{MN^\T} = ML_\rN^\T N_\perp^\T + M(N^\T N)K(N^\T N)^{-1}N^\T + M_\perp L_\rM N^\T +
M(M^\T M)^{-1} K (M^\T M) N^\T.
\end{equation}
\begin{subequations}  \label{eq:LK}
Multiplying~\eqref{eq:dX=LK} on the left by $(M^\T M)^{-1}M^\T$ yields
\begin{equation}  \label{eq:LN}
L_\rN^\T = (M^\T M)^{-1}M^\T\dX_{MN^\T} N_\perp,
\end{equation}
multiplying~\eqref{eq:dX=LK} on the right by $N(N^\T N)^{-1}$ yields
\begin{equation}  \label{eq:LM}
L_\rM = M_\perp^\T\dX_{MN^\T} N(N^\T N)^{-1},
\end{equation}
and multiplying~\eqref{eq:dX=LK} on the left by $M^\T$ and on the right by
$N$ yields
\begin{equation}  \label{eq:K}
M^\T\dX_{MN^\T} N = M^\T MN^\T NK + KM^\T MN^\T N.
\end{equation}
\end{subequations}
Replacing~\eqref{eq:LK} into~\eqref{eq:dXMN=LK} yields
\begin{subequations}  \label{eq:dXMN-1}
\begin{gather}
\dX_{\rM(M,N)} = M_\perp M_\perp^\T\dX_{MN^\T} N(N^\T N)^{-1} + M(M^\T M)^{-1}KM^\T M
\\ \dX_{\rN(M,N)} = N_\perp N_\perp^\T \dX_{MN^\T}^\T M(M^\T M)^{-1} +
N(N^\T N)^{-1}K^\T N^\T N.
\end{gather}
\end{subequations}
We can further exploit the identities $M_\perp M_\perp^\T =
I-M(M^\T M)^{-1}M^\T$, and likewise for $N$, to rewrite~\eqref{eq:dXMN-1}
as
\begin{subequations}  \label{eq:dXMN-2}
\begin{gather}
\dX_{\rM(M,N)} = (\dX_{MN^\T} N - MN^\T NK)(N^\T N)^{-1}
\\ \dX_{\rN(M,N)} = (\dX_{MN^\T}^\T M - NM^\T M K^\T)(M^\T M)^{-1}.
\end{gather}
\end{subequations}
This result is formalized as follows:
\begin{prpstn}   \label{prp:H-lift}
Consider the submersion $\pi$~\eqref{eq:pi} and the
horizontal distribution~\eqref{eq:H-def}. Let $(M,N)\in\Rp$ and let $\dX_{MN^\T}$ be in
$\T_{MN^\T}\Mp$. Then the horizontal lift of $\dX_{MN^\T}$ at $(M,N)$ is $\dX_{(M,N)}
= (\dX_{\rM(M,N)},\dX_{\rN(M,N)})$ given by~\eqref{eq:dXMN-2}, where $K$ is the
solution of the Sylvester equation~\eqref{eq:K}. 
\end{prpstn}

\subsection{Constitutive equation of horizontal lifts}

A horizontal lift $\dX_{(M,N)}$ fully specifies $\dX_{MN^\T}
= \D\pi(M,N)[\dX_{(M,N)}] \in \T_{MN^\T}\Mp$
as well as its horizontal lift at any other point of the
fiber $\pi^{-1}(MN^\T)$~\eqref{eq:pi-1}. Let us obtain an expression for
$\dX_{(MR,NR^{-\T})}$ in terms of $\dX_{(M,N)}$. The
expression~\eqref{eq:dXMN-2} of horizontal lifts yields after routine
manipulations 
\begin{equation}  \label{eq:lift-eq}
\dX_{\rM(MR,NR^{-\T})} = \dX_{\rM(M,N)} R, \quad \dX_{\rN(MR,NR^{-\T})} =
\dX_{\rN(M,N)} R^{-\T}.
\end{equation}
We have obtained:
\begin{prpstn}  \label{prp:lift-eq}
Consider the submersion $\pi$~\eqref{eq:pi} and the
horizontal distribution~\eqref{eq:H-def}. Then a vector field
$\Rp\ni(M,N)\mapsto \dX_{(M,N)}\in\TRp$ is a horizontal lift if and
only if~\eqref{eq:lift-eq} holds for all $(M,N)\in\Rp$ and all
$R\in\GL(p)$. 
\end{prpstn}

\subsection{Riemannian submersion}

Routine manipulations using~\eqref{eq:lift-eq} yield that
$\bar{g}$~\eqref{eq:bar-g} satisfies the invariance
condition~\eqref{eq:g=g}. Hence there is a (unique) Riemannian metric $g$ on $\Mp$ that makes 
\begin{equation}  \label{eq:pi-RS}
\pi:(\Rp,\bar{g})\to(\Mp,g): (M,N)\mapsto MN^\T
\end{equation}
a Riemannian submersion. 
The Riemannian metric
$g$ is consistently defined by
\begin{equation}  \label{eq:g}
g_{MN^\T}(\dX_{MN^\T},\cX_{MN^\T}) := \bar{g}_{(M,N)}(\dX_{(M,N)},\cX_{(M,N)}).
\end{equation}

\subsection{Horizontal projection}
\label{sec:Ph}

We will need an expression for the projection $P^\rh_{(M,N)}(\dM,\dN)$
of $(\dM,\dN)\in\TRp$
onto the horizontal space~\eqref{eq:H-def} along the vertical
space~\eqref{eq:V}. 

Since the projection is along the vertical space, we have
\begin{equation}  \label{eq:Ph}
P^\rh_{(M,N)}(\dM,\dN) = (\dM+M\dR,\dN-N\dR^\T)
\end{equation}
for some
$\dR\in\mathbb{R}^{p\times p}$. It remains to obtain $\dR$ by imposing
horizontality of~\eqref{eq:Ph}. Since horizontal vectors are
characterized by~\eqref{eq:H-cond}, we find that~\eqref{eq:Ph} is
horizontal if and only if
\[
M^\T(\dM+M\dR)(M^\T M)^{-1} = (N^\T N)^{-1}(\dN^\T-\dR N^\T)N,
\]
that is,
\[
M^\T M\dR(M^\T M)^{-1} + (N^\T N)^{-1}\dR N^\T N = -M^\T\dM(M^\T M)^{-1} +
(N^\T N)^{-1}\dN^\T N,
\]
which can be rewritten as the Sylvester equation
\begin{equation}  \label{eq:dR-pj}
N^\T NM^\T M\dR + \dR N^\T NM^\T M = -N^\T NM^\T\dM +
\dN^\T NM^\T M.
\end{equation}
In summary:
\begin{prpstn}  \label{prp:Ph}
The projection $P^\rh_{(M,N)}(\dM,\dN)$ of $(\dM,\dN)\in\TRp$
onto the horizontal space~\eqref{eq:H-def} along the vertical
space~\eqref{eq:V} is given by~\eqref{eq:Ph} where $\dR$ is
the solution of the Sylvester equation~\eqref{eq:dR-pj}.
\end{prpstn}

\subsection{Riemannian connection on the total space}
\label{sec:bar-nabla}

Since the chosen Riemannian metric $\bar{g}$~\eqref{eq:bar-g} on the total space
$\Rp$ is not the Euclidean metric~\eqref{eq:Euc-bar-g}, it can be
expected that the Riemannian connection on $(\Rp,\bar{g})$ is not
the plain differential. We show that this is indeed the case and we
provide a formula for the Riemannian connection $\bar\nabla$ on $(\Rp,\bar{g})$. The
motivation for obtaining this formula is that the Riemannian Newton equation on $(\Mp,g)$
requires the Riemannian connection on $(\Mp,g)$, which is readily
obtained from $\bar\nabla$ as we will see in
Section~\ref{sec:nabla}. The general theory of Riemannian connections
(also called Levi-Civita connections) can be found in~\cite[\S
  5.3]{AMS2008} or in any Riemannian geometry textbook such
as~\cite{dC92}.

The development relies on Koszul's formula 
\begin{equation}  \label{eq:Koszul}
2g(\nabla_\chi\eta,\xi) = \partial_\chi g(\eta,\xi) + \partial_\eta g(\chi,\xi) -
\partial_\xi g(\chi,\eta) + g([\chi,\eta],\xi) - g([\chi,\xi],\eta) - g([\eta,\xi],\chi).
\end{equation}

After lengthy but routine calculations, we obtain the following expression for the
Riemannian connection $\bar\nabla$ on $(\Rp,\bar{g})$:
\begin{subequations}  \label{eq:bar-nabla}
\begin{multline}
\left(\bar\nabla_\dX \dY\right)_\rM 
= \partial_\dX \dY_\rM  - \dY_\rM (M^\T M)^{-1}
\sym(\dX_\rM^\T M) - \dX_\rM (M^\T M)^{-1}\sym(\dY_\rM^\T M) \\ + M(M^\T M)^{-1}
\sym(\dX_\rM^\T\dY_\rM)
\end{multline}
and
\begin{multline}
\left(\bar\nabla_\dX \dY\right)_\rN = \partial_\dX \dY_\rN - \dY_\rN (N^\T N)^{-1}
\sym(\dX_\rN^\T N) - \dX_\rN (N^\T N)^{-1}\sym(\dY_\rN^\T N) \\ + N(N^\T N)^{-1}
\sym(\dX_\rN^\T\dY_\rN),
\end{multline}
\end{subequations}
for all $(M,N)\in\Rp$, all $\dX\in \T_{(M,N)}\Rp$ and all tangent vector fields $\dY$ on $\Rp$.

\subsection{Connection on the quotient space}
\label{sec:nabla}

Let $\nabla$ denote the Riemannian connection on the quotient space $\Mp$ endowed with the Riemannian metric $g$~\eqref{eq:g}.
A classical result in the theory of Riemannian submersions (see~\cite[Lemma~1]{ONe1966} or~\cite[\S 5.3.4]{AMS2008}) states that 
\[
(\nabla_{\dX_{MN^\T}}\dY)_{(M,N)} = P^\rh_{(M,N)} (\bar\nabla_{\dX_{(M,N)}}\dY),
\]
for all $\dX_{MN^\T}\in\T_{MN^\T}\Mp$ and all tangent tangent vector fields $\dY$ on $\Mp$.
That is, the horizontal lift of the Riemannian connection of the
quotient space is given by the horizontal projection~\eqref{eq:Ph} of
the Riemannian connection~\eqref{eq:bar-nabla} of the total
space. (The tangent vector field $Y$ on the right-hand side denotes the
horizontal lift of the tangent vector field $Y$ of the left-hand side.) 

\subsection{Riemannian Newton equation}

For a real-valued function $f$ on a Riemannian manifold $\mathcal{M}$ with Riemannian metric $g$, we let $\grad\,f(x)$ denote the gradient of $f$ at $x\in\mathcal{M}$---defined
as the unique tangent vector to $\mathcal{M}$ at $x$ that satisfies
$g_x(\grad\,f(x),\xi_x)=\mathrm{D}f(x)[\xi_x]$ for all $\xi_x\in
\T_x\mathcal{M}$---and the plain Riemannian Newton equation is given by 
\[
\nabla_{\eta_x} \grad\,f = -\grad\,f(x)
\] 
for the unknown $\eta_x\in \T_x\mathcal{M}$, where $\nabla$ stands for the Riemannian connection; see, e.g.,~\cite[\S 6.2]{AMS2008}.

We now turn to the manifold $\Mp$ endowed with the Riemannian metric
$g$~\eqref{eq:g} and we obtain an expression of the Riemannian Newton equation by means of its horizontal lift through the Riemannian submersion $\pi$~\eqref{eq:pi-RS}. 
First, on the total space $\Rp$ endowed with the Riemannian metric $\bar{g}$~\eqref{eq:bar-g}, 
we readily obtain the following expression for the gradient of $\bar{f}$~\eqref{eq:bar-f}:
\[
\grad\,\bar{f}(M,N) = (\partial_\rM \bar{f}(M,N) M^\T M, \partial_\rN
\bar{f}(M,N) N^\T N),
\]
where $\partial_\rM \bar{f}(M,N)$ denotes the Euclidean (i.e., classical) gradient of
$\bar{f}$ with respect to its first argument, i.e., $(\partial_\rM
\bar{f}(M,N))_{i,j} =
\frac{\mathrm{d}}{\mathrm{d}t}\bar{f}(M+te_i^{}e_j^\T,N)|_{t=0}$, and
likewise for $\partial_\rN \bar{f}(M,N)$ with the second argument. Then the horizontal lift of the Newton
equation at a point $(M,N)$ of the total space $\Rp$, for 
the unknown $X_{(M,N)}$ in the horizontal space
$\mathcal{H}_{(M,N)}$~\eqref{eq:H-def}, is
\begin{equation}  \label{eq:N}
P^\rh_{(M,N)} (\bar\nabla_{X_{(M,N)}}\grad\,\bar{f}) =
-\grad\,\bar{f}(M,N),
\end{equation}
where $P^\rh$ is the horizontal projection given in
Section~\ref{sec:Ph} and $\bar\nabla$ is the Riemannian connection on
$(\Rp,\bar{g})$ given in Section~\ref{sec:bar-nabla}. To obtain~\eqref{eq:N}, we
have used the fact (see~\cite[(3.39)]{AMS2008}) that $\grad\,f(M,N) =
\grad\,\bar{f}(M,N)$, where the left-hand side denotes the horizontal
lift of $\grad\,f(MN^\T)$ at $(M,N)$.

Intimidating as it may be in view of the expressions of $P^\rh$ and
$\bar\nabla$, the Newton equation~\eqref{eq:N} is nevertheless merely
a linear system of equations. Indeed, $X_{(M,N)} \mapsto P^\rh_{(M,N)}
(\bar\nabla_{X_{(M,N)}}\grad\,\bar{f})$ is a linear transformation of
the horizontal space $\mathcal{H}_{(M,N)}$. Thus~\eqref{eq:N} can be solved using
``matrix-free'' linear solvers such as GMRES. Moreover, in addition to computing the Euclidean gradient of $\bar{f}$ and the Euclidean derivative of the Euclidean gradient of $f$ along $X_{(M,N)}$, computing $P^\rh_{(M,N)} (\bar\nabla_{X_{(M,N)}}\grad\,\bar{f})$ requires only
$O(p^2(m+n+p))$ flops. 

\subsection{Newton's method}

In order to spell out on $(\Mp,g)$ the Riemannian Newton method as defined in~\cite[\S 6.2]{AMS2008}, the
last missing ingredient is a retraction $R$ that turns the Newton
vector $\dX_{MN^\T}$ into an updated iterate $R_{MN^\T}\dX_{MN^\T}$ in $\Mp$. The general definition of a retraction can be found
in~\cite[\S 4.1]{AMS2008}. 

The quintessential retraction on a Riemannian manifold is the Riemannian exponential; see~\cite[\S
  5.4]{AMS2008}. However, computing the Riemannian exponential amounts to
solving the differential equation $\nabla_{\dX}\dX=0$, which may
not admit a closed-form solution. In the case of $(\Mp,g)$, we are not
aware of such a closed-form solution, and this makes the exponential
retraction impractical.

Fortunately, other retractions are readily available. A retraction on
$\Mp$ is given by
\begin{equation}  \label{eq:R}
R_{MN^\T}(\dX_{MN^\T}) := (M+\dX_{\rM(M,N)})(N+\dX_{\rN(M,N)})^\T,
\end{equation}
where $\dX_{\rM(M,N)}$ and $\dX_{\rN(M,N)}$ are horizontal lifts as defined in
Proposition~\ref{prp:H-lift}. It is readily checked that the
definition is consistent, i.e., it depends on $MN^\T$ and not on the
specific choices of $(M,N)$ in the fiber~\eqref{eq:pi-1}.

With all these elements in place, we can describe Newton's method as
follows.
\begin{thrm}[Riemannian Newton on $\Mp$ with Riemannian metric~\eqref{eq:g}]  \label{thm:N}
Let $f$ be a real-valued function on the Riemannian manifold
$\Mp$~\eqref{eq:Mp}, endowed with the Riemannian metric
$g$~\eqref{eq:g}, with the associated Riemannian connection, and with
the retraction~\eqref{eq:R}. Then the Riemannian Newton method for $f$
maps $MN^\T\in\Mp$ to $(M+\dX_\rM)(N+\dX_\rN)^\T$, where
$(\dX_\rM,\dX_\rN)$ is the solution $\dX_{(M,N)}$ of the Newton
equation~\eqref{eq:N}.
\end{thrm} 

Note that, in practice, it is not necessary to form $MN^\T$. Given an
initial point $M_0^{}N_0^\T$, one can instead generate a sequence
$\{(M_k,N_k)\}$ in $\Rp$ by applying the iteration map $(M,N)\mapsto
(M+\dX_\rM,N+\dX_\rN)$. The Newton sequence on $\Mp$ is then
$\{M_k^{}N_k^\T\}$, and it depends on $M_0^{}N_0^\T$ but not on the particular
$M_0$ and $N_0$.

The following convergence result follows directly from the general convergence analysis of the Riemannian Newton method~\cite[Theorem~6.3.2]{AMS2008}. A critical point of $f:\Mp\to\mathbb{R}$ is a point $X_*$ where $\grad\,f(X_*)=0$. It is termed \emph{nondegenerate} if the Hessian $\T_{X_*}\Mp\ni\dX\mapsto \nabla_\dX\grad\,f\in \T_{X_*}\Mp$ is invertible. These definitions do not depend on the Riemannian metric nor on the affine connection $\nabla$.
\begin{thrm}[quadratic convergence]  \label{thm:quad}
Let $X_*$ be a nondegenerate critical point of $f$. Then there exists a neighborhood $\mathcal{U}$ of $X_*$ in $\Mp$ such that, for all initial iterate $X_0\in\mathcal{U}$, the iteration described in Theorem~\ref{thm:N} generates an infinite sequence $\{X_k\}$ converging superlinearly (at least quadratically) to $X_*$. 
\end{thrm}

\section{$\Mp$ as a Riemannian quotient manifold with an orthonormal factor}
\label{sec:R-o}  

We now follow the second plan of action mentioned at the end of
Section~\ref{sec:intro}. Bear in mind that the meaning of much of the
notation introduced above will be superseded by new definitions below.

\subsection{A smaller total space}

Let 
\begin{equation}  \label{eq:St}
\St{p}{m} = \{M\in\mathbb{R}^{m\times p}: M^\T M=I_p\},
\end{equation}
denote the \emph{Stiefel manifold} of orthonormal
$m\times p$ matrices.
For all $X\in\Mp$, there exists $(M,N)$ with $M$ orthonormal such that
$X=MN^\T$. To see this, take $(M,N)\in\Rp$ such that $X=MN^\T$, let
$M=QR$ be a QR decomposition of $M$, where $R$ is invertible since $M$
has full rank, and observe that $X=MR^{-1}(NR^\T)^\T = Q(NR^\T)^\T$. Hence
\begin{equation}  \label{eq:R-o:pi}
\pi: \SRp \to \Mp: (M,N)\mapsto MN^\T
\end{equation}
is a smooth surjective map between two manifolds. 

As in Section~\ref{sec:quot}, but now with the restricted total space $\SRp$, we show that $\pi$~\eqref{eq:R-o:pi} is a
submersion.  The tangent space at $M$ to $\St{p}{m}$ is given by
(see~\cite[Example~3.5.2]{AMS2008})
\begin{align*}
\T_M\St{p}{m} &= \{\dM\in\mathbb{R}^{m\times p}: M^\T\dM + \dM^\T M = 0\}
\\ &= \{M\Omega + M_\perp W: \Omega = -\Omega^\T\in\mathbb{R}^{p\times p},
W\in\mathbb{R}^{(m-p)\times p}\},
\end{align*}
and we have
\[
\T_{(M,N)}\SRp = (\T_M\St{p}{m}) \times \mathbb{R}^{n\times p}.
\] 
For all $(M,N)\in\SRp$ and all $(\dot{M},\dot{N})\in \T_{(M,N)}\SRp$,
we have that
$\mathrm{D}\pi(M,N)[(\dot{M},\dot{N})] = \dot{M}N^\T +
M\dot{N}^\T$. 
Here again, we can work in a coordinate system where $M
= \begin{bmatrix} I & 0 \end{bmatrix}^\T$ and $N = \begin{bmatrix} I &
  0 \end{bmatrix}^\T$. We have that
$\{\mathrm{D}\pi(M,N)[(\dot{M},\dot{N})]: (\dot{M},\dot{N})\in
\T_{(M,N)}\SRp\} = \{\left[\begin{smallmatrix}\Omega+N_1^\T & N_2^\T \\ W
    & 0\end{smallmatrix}\right]:
\Omega=-\Omega^\T\in\mathbb{R}^{p\times p}, N_1\in\mathbb{R}^{p\times
  p}, N_2\in\mathbb{R}^{(n-p)\times p}, W\in\mathbb{R}^{(m-p)\times
  p}\}$, a linear subspace of dimension $p^2+(n-p)p+(m-p)p =
p(m+n-p)$, which is the dimension of $\Mp$. Hence
$\pi$~\eqref{eq:R-o:pi} is a submersion.

The fiber of $\pi$~\eqref{eq:R-o:pi} at $MN^\T$ is now
\begin{equation}  \label{eq:R-o:pi-1}
\pi^{-1}(MN^\T) = \{(MR,NR): R\in\mathrm{O}(p)\},
\end{equation}
where 
\[
\mathrm{O}(p) = \{R\in\mathbb{R}^{p\times p}:R^\T R=I_p\}
\]
denotes the orthogonal group of degree $p$. 

The vertical space $\mathcal{V}_{(M,N)}$ at a point $(M,N)\in\Rp$,
i.e., the tangent space to
the fiber $\pi^{-1}(MN^\T)$ at $(M,N)$, is given by
\begin{equation} \label{eq:R-o:V}
\mathcal{V}_{(M,N)} = \{(M\Omega,N\Omega): \Omega = -\Omega^\T\in\mathbb{R}^{p\times p}\}.
\end{equation}

\subsection{Riemannian metric in total space}

We consider $\SRp$ as a Riemannian submanifold of the Euclidean space $\TRp$. This endows
$\SRp$ with the Riemannian metric $\bar{g}$ defined by
\begin{equation}  \label{eq:R-o:bar-g}
\bar{g}_{(M,N)}\left( (\dM,\dN), (\cM,\cN) \right) := \trace
\left(\dM^\T\cM + \dN^\T\cN \right)
\end{equation}
for all $(\dM,\dN)$ and $(\cM,\cN)$ in $\T_{(M,N)}\SRp$. 

Adapting the rationale of Section~\ref{sec:R-no}, we will obtain in Section~\ref{sec:R-o:RS} below that, with this $\bar{g}$, $\pi$~\eqref{eq:R-o:pi} can be turned into a Riemannian submersion. 

\subsection{Horizontal space}

The horizontal space $\mathcal{H}_{(M,N)}$ is the orthogonal
complement to $\mathcal{V}_{(M,N)}$~\eqref{eq:R-o:V} in
$\T_{(M,N)}\SRp$ with respect to $\bar{g}$~\eqref{eq:R-o:bar-g}. The
following propositions are equivalent:
\begin{gather}
(\dot{M},\dot{N})\in\mathcal{H}_{(M,N)}, \nonumber
\\ \dM \in \T_M\St{p}{n},\ \dN \in \mathbb{R}^{n\times p},\
\mathrm{tr}(\tilde\Omega^\T(M^\T\dM+N^\T\dN))=0, \forall \tilde\Omega =
-\tilde\Omega^\T,  \nonumber
\\ M^\T\dM = -(M^\T\dM)^\T,\ M^\T\dM + N^\T\dN = (M^\T\dM + N^\T\dN)^\T,  \label{eq:R-o:H-ss}
\\ \dM = M\Omega + M_\perp W,\
\dN = N(N^\T N)^{-1}(-\Omega+S)+N_\perp L, \nonumber
\end{gather}
with $W\in\mathbb{R}^{(m-p)\times p}$, $\Omega =
-\Omega^\T\in\mathbb{R}^{p\times p}$, $S=S^\T\in\mathbb{R}^{p\times
  p}$, $L\in\mathbb{R}^{(n-p)\times p}$. In summary,
\begin{equation}  \label{eq:R-o:H-def}
\mathcal{H}_{(M,N)} = \{(\dM,\dN): M^\T\dM = -(M^\T\dM)^\T,\ M^\T\dM + N^\T\dN = (M^\T\dM + N^\T\dN)^\T\}.
\end{equation}

\subsection{Horizontal lift}
\label{sec:R-o:hor-lift}

Proceeding as in Section~\ref{sec:R-no-hor-lift} but now with the
horizontal space~\eqref{eq:R-o:H-def} and taking into account that
$M^\T M=I$, we obtain that the horizontal lift of
$\dX_{MN^\T}\in\T_{MN^\T}\SRp$ is given by
\begin{subequations}  \label{eq:R-o:hor-lift-1}
\begin{gather}
\dX_{\rM(M,N)} = M\Omega + M_\perp M_\perp^\T\dX_{MN^\T}N(N^\T N)^{-1}
\\ \dX_{\rN(M,N)} = N(N^\T N)^{-1}(S-\Omega) + N_\perp N_\perp^\T\dX_{MN^\T}^\T M
\\ \text{where} \quad \Omega(N^\T N+I)+S = M^\T\dX_{MN^\T}N,\ \Omega=-\Omega^\T,\ S=S^\T.  \label{eq:R-o:OS1}
\end{gather}
\end{subequations}
Equation~\eqref{eq:R-o:OS1} is equivalent to
\begin{subequations}   \label{eq:R-o:OS2}
\begin{gather}
\Omega (N^\T N+I) + (N^\T N+I)\Omega = M^\T\dX_{MN^\T}N - N^\T\dX_{MN^\T}^\T M, \label{eq:R-o:OS2O}
\\ S = M\dX_{MN^\T}N - \Omega(N^\T N+I).  \label{eq:R-o:OS2S}
\end{gather}
\end{subequations}
As for the first two equations of~\eqref{eq:R-o:hor-lift-1},
using~\eqref{eq:R-o:OS1}, they can
be rewritten as 
\begin{subequations}  \label{eq:R-o:hor-lift-2}
\begin{gather}
\dX_{\rM(M,N)} = \dX_{MN^\T}N(N^\T N)^{-1} - M(\Omega+S)(N^\T N)^{-1}
\\ \dX_{\rN(M,N)} = \dX_{MN^\T}^\T M + N\Omega.
\end{gather}
\end{subequations}
In summary,
\begin{prpstn}   \label{prp:R-o:H-lift}
Consider the submersion $\pi$~\eqref{eq:R-o:pi} and the
horizontal distribution~\eqref{eq:R-o:H-def}. Let $(M,N)\in\SRp$ and let $\dX_{MN^\T}\in\T_{MN^\T}\Mp$. Then the horizontal lift of $\dX_{MN^\T}$ at $(M,N)$ is $\dX_{(M,N)}
= (\dX_{\rM(M,N)},\dX_{\rN(M,N)})$ given by~\eqref{eq:R-o:hor-lift-2}, where $\Omega$ is the
solution of the Sylvester equation~\eqref{eq:R-o:OS2O} and $S$ is
given by~\eqref{eq:R-o:OS2S}. 
\end{prpstn}

\subsection{Constitutive equation of horizontal lifts}

From Proposition~\ref{prp:R-o:H-lift}, routine manipulations lead to the
following constitutive equation for horizontal lifts:
\begin{equation}  \label{eq:R-o:lift-eq}
\dX_{\rM(MR,NR)} = \dX_{\rM(M,N)} R, \quad \dX_{\rN(MR,NR)} =
\dX_{\rN(M,N)} R.
\end{equation}
Hence we have the following counterpart of Proposition~\ref{prp:lift-eq}.
\begin{prpstn}  \label{prp:R-o:lift-eq}
Consider the submersion $\pi$~\eqref{eq:R-o:pi} and the
horizontal distribution~\eqref{eq:R-o:H-def}. Then a tangent vector field
$\SRp\ni(M,N)\mapsto \dX_{(M,N)}\in \T_{(M,N)}\SRp$ is a horizontal lift if and
only if~\eqref{eq:R-o:lift-eq} holds for all $(M,N)\in\SRp$ and all
$R\in\mathrm{O}(p)$. 
\end{prpstn}

\subsection{Riemannian submersion}
\label{sec:R-o:RS}

From Proposition~\ref{prp:R-o:lift-eq} and the properties of the
trace, it is direct that $\bar{g}$~\eqref{eq:R-o:bar-g} satisfies the
invariance condition
\begin{equation}  \label{eq:R-o:g=g}
\bar{g}_{(M,N)}(\dX_{(M,N)},\cX_{(M,N)}) =
\bar{g}_{(MR,NR)}(\dX_{(MR,NR)},\cX_{(MR,NR)}).
\end{equation}
Hence one consistently defines a Riemannian metric $g$ on $\Mp$ by
\begin{equation}  \label{eq:R-o:g}
g_{MN^\T}(\dX_{MN^\T},\cX_{MN^\T}) =
\bar{g}_{(M,N)}(\dX_{(M,N)},\cX_{(M,N)}),
\end{equation}
and $\pi:(\SRp,\bar{g})\to(\Mp,g)$ is a Riemannian submersion.

\subsection{Horizontal projection}
\label{sec:R-o:Ph}

We now obtain an expression for the projection $P^\rh_{(M,N)}(\dM,\dN)$
of $(\dM,\dN)\in \T_{(M,N)}\SRp$
onto the horizontal space~\eqref{eq:R-o:H-def} along the vertical
space~\eqref{eq:R-o:V}. 
Since the projection is along the vertical space, we have
\begin{equation}  \label{eq:R-o:Ph}
P^\rh_{(M,N)}(\dM,\dN) = (\dM+M\Omega,\dN+N\Omega)
\end{equation}
for some
$\Omega=-\Omega^\T\in\mathbb{R}^{p\times p}$. It remains to obtain $\Omega$ by imposing
horizontality of~\eqref{eq:R-o:Ph}. The characterization of horizontal
vectors given in~\eqref{eq:R-o:H-ss} yields the Sylvester equation
\begin{equation}  \label{eq:R-o:Om-pj}
(N^\T N+I)\Omega + \Omega (N^\T N+I) = \dM^\T M - M^\T\dM + \dN^\T N - N^\T\dN.
\end{equation}
In summary:
\begin{prpstn}  \label{prp:R-o:Ph}
The projection $P^\rh_{(M,N)}(\dM,\dN)$ of $(\dM,\dN)\in \T_{(M,N)}\SRp$
onto the horizontal space~\eqref{eq:R-o:H-def} along the vertical
space~\eqref{eq:R-o:V} is given by~\eqref{eq:R-o:Ph} where $\Omega$ is
the solution of the Sylvester equation~\eqref{eq:R-o:Om-pj}.
\end{prpstn}

\subsection{Riemannian connection on the total space}
\label{sec:R-o:bar-nabla}

Let $P^{\StX}_M$ denote the orthogonal projection from
$\mathbb{R}^{m\times p}$ onto
$\T_M\St{p}{n}$, given by (see~\cite[Example~5.3.2]{AMS2008})
\begin{equation}  \label{eq:R-o:PSt}
P^{\StX}_M\dM = (I-MM^\T)\dM + M\skew(M^\T\dM) = \dM - M\sym(M^\T\dM),
\end{equation}
where $\skew(Z):=\frac12(Z-Z^\T)$ and $\sym(Z):=\frac12(Z+Z^\T)$. We
also let $P^{\StX\times\mathbb{R}}_{(M,N)}$ denote the orthogonal projection
from $\TRp$ onto $\T_{(M,N)}\SRp$, given by
\begin{equation}  \label{eq:R-o:PStR}
P^{\StX\times\mathbb{R}}_{(M,N)}(\dM,\dN) = (P^{\StX}_M\dM, \dN).
\end{equation}

Since $\SRp$, endowed with the Riemannian metric
$\bar{g}$~\eqref{eq:R-o:bar-g}, is a Riemannian submanifold of the
Euclidean space $\TRp$, a classical result of Riemannian geometry
(see~\cite[\S 5.3.3]{AMS2008}) yields that the Riemannian connection
$\bar\nabla$ on $(\SRp,\bar{g})$ is given by
\[
\bar\nabla_\dX \dY = P^{\StX\times\mathbb{R}}_{(M,N)} \partial_\dX \dY,
\]
that is,
\begin{subequations}  \label{eq:R-o:bar-nabla}
\begin{align}
\left(\bar\nabla_\dX \dY\right)_\rM 
&= P^{\StX}_M(\partial_\dX \dY_\rM)
\\ \left(\bar\nabla_\dX \dY\right)_\rN &= \partial_\dX \dY_\rN 
\end{align}
\end{subequations}
for all $(M,N)\in\SRp$, all $\dX\in \T_{(M,N)}\SRp$ and all vector
fields $\dY$ on $\SRp$.  

\subsection{Connection on the quotient space}
\label{sec:R-o:nabla}

As in Section~\ref{sec:nabla}, we can now provide an expression for
the Riemannian connection $\nabla$ on the manifold $\Mp$ endowed with
the Riemannian metric $g$~\eqref{eq:R-o:g}:
\begin{align*}
(\nabla_{\dX_{MN^\T}}\dY)_{(M,N)} &= P^\rh_{(M,N)}
  (\bar\nabla_{\dX_{(M,N)}}\dY)
\\ &= P^\rh_{(M,N)} P^{\StX\times\mathbb{R}}_{(M,N)} \partial_\dX \dY,
\end{align*}
with $P^\rh$ as in~\eqref{eq:R-o:Ph} and $P^{\StX\times\mathbb{R}}$ as
in~\eqref{eq:R-o:PStR}. (Observe that $\dY$ of the right-hand side is the horizontal lift of $\dY$ of the left-hand side.)  

\subsection{Riemannian Newton equation}

Given $f:\Mp\to\mathbb{R}$, define
$\bar{f}=f\circ\pi$, i.e.,
\[
\bar{f}: \SRp\to\mathbb{R}: (M,N)\mapsto f(MN^\T),
\]
and define
\[
\bar{\bar{f}}: \Rp\to\mathbb{R}: (M,N)\mapsto f(MN^\T).
\]
Let $\grad\,\bar{\bar{f}}$ denote the Euclidean gradient of
$\bar{\bar{f}}$. We have (see~\cite[(3.37)]{AMS2008}) 
\begin{equation}  \label{eq:R-o:gbf}
\grad\,\bar{f}(M,N) = P^{\StX\times\mathbb{R}}_{(M,N)} \grad\,\bar{\bar{f}}(M,N)
\end{equation}
and (see~\cite[(3.39)]{AMS2008})
\[
\grad\,f(M,N) = \grad\,\bar{f}(M,N),
\]
where the left-hand side stands for the horizontal lift at $(M,N)$ of
$\grad\,f(MN^\T)$.

We can now obtain the counterpart of the (lifted) Newton equation~\eqref{eq:N} with normalization
on the $\rM$ factor: 
\begin{equation}  \label{eq:R-o:N}
P^\rh_{(M,N)} (\bar\nabla_{X_{(M,N)}}\grad\,\bar{f}) =
-\grad\,\bar{f}(M,N),
\end{equation}
where $P^\rh$ is the horizontal projection given in
Section~\ref{sec:R-o:Ph}, $\bar\nabla$ is the Riemannian connection on
$(\Rp,\bar{g})$ given in Section~\ref{sec:R-o:bar-nabla}, and
$\grad\,\bar{f}$ is obtained from the Euclidean gradient of
$\bar{\bar{f}}$ from~\eqref{eq:R-o:gbf}.

The Newton equation~\eqref{eq:R-o:N} can be considered less intricate than in
the non-orthonormal case~\eqref{eq:N} because the expression for
$\bar\nabla$ in~\eqref{eq:R-o:bar-nabla} is simpler than
in~\eqref{eq:bar-nabla}. In any case, the discussion that
follows~\eqref{eq:N} applies equally: the Newton equation is merely a
linear system of equations, and the Riemannian overhead requires only
$O(p^2(m+n+p))$ flops. 

\subsection{Newton's method}

Another reward that comes with the orthonormalization of the M factor
is that the Riemannian exponential with respect to $g$~\eqref{eq:R-o:g} admits a
closed-form expression. First, we point out that, in view
of~\cite[\S 2.2.2]{EAS98}, the Riemannian
exponential on $\SRp$ for the Riemannian metric $\bar{g}$~\eqref{eq:R-o:bar-g}
is given by
\begin{equation}  \label{eq:R-o:SRp-Exp}
\Exp_{(M,N)}(\dM,\dN) = ( \begin{bmatrix}M & \dM\end{bmatrix}
  \exp\begin{bmatrix}A & -S \\ I & A\end{bmatrix} I_{2p,p}
  \exp(-A),\ N+\dN ),
\end{equation}
where $A:=M^\T\dM$ and $S:=\dM^\T\dM$, and where $\exp$ stands for the matrix
exponential (\texttt{expm} in Matlab). Second, since by~\cite[Corollary~7.46]{ONe1983} horizontal
geodesics in $(\SRp,\bar{g})$ map to geodesics in $(\Mp,g)$, we have
that
\begin{equation}  \label{eq:R-o:Mp-Exp}
\Exp_{MN^\T}(\dX_{MN^\T}) = \pi(\Exp_{(M,N)}(\dX_{\rM(M,N)},\dX_{\rN(M,N)})), 
\end{equation}
with $(\dX_{\rM(M,N)},\dX_{\rN(M,N)})$ as in
Proposition~\ref{prp:R-o:H-lift}. (In~\eqref{eq:R-o:Mp-Exp}, $\Exp$ on the right-hand side is given by~\eqref{eq:R-o:SRp-Exp} and $\Exp$ on the left-hand side denotes the Riemannian exponential of $(\Mp,g)$.)

Observe that the matrix exponential is applied in~\eqref{eq:R-o:SRp-Exp} to matrices of size $2p\times 2p$ and $p\times p$; hence, when $p\ll m$, the cost of computing the M component of~\eqref{eq:R-o:SRp-Exp} is comparable to the cost of computing the simple sum $M+\dM$. Note also that, in practice, the M component of the Newton iterates may gradually depart from orthonormality due to the accumulation of numerical errors; a remedy is to restore orthonormality by taking the Q factor of the unique QR decomposition where the diagonal of the R factor is positive.

We can now formally describe Newton's method in the context of this
Section~\ref{sec:R-o}.
\begin{thrm}[Riemannian Newton on $\Mp$ with Riemannian metric~\eqref{eq:R-o:g}]  \label{thm:R-o:N}
Let $f$ be a real-valued function on the Riemannian manifold
$\Mp$~\eqref{eq:Mp}, endowed with the Riemannian metric
$g$~\eqref{eq:R-o:g}, with the associated Riemannian connection, and
with the exponential retraction~\eqref{eq:R-o:Mp-Exp}. Then the
Riemannian Newton method for $f$ maps $MN^\T\in\Mp$ to
$\pi(\Exp_{(M,N)}(\dX_\rM,\dX_\rN))$, where $\pi$ is given
in~\eqref{eq:R-o:pi}, $\Exp$ is defined
in~\eqref{eq:R-o:SRp-Exp}, and $(\dX_\rM,\dX_\rN)$ is the solution
$\dX_{(M,N)}$ of the Newton equation~\eqref{eq:R-o:N}.
\end{thrm} 

The quadratic convergence result in Theorem~\ref{thm:quad} still
holds, replacing the reference to Theorem~\ref{thm:N} by a reference
to Theorem~\ref{thm:R-o:N}. 

\section{Conclusion}
\label{sec:conc}

We have reached the end of a technical hike that led us to give in
Theorem~\ref{thm:R-o:N} what is, to the best of our knowledge, the
first closed-form description of a purely Riemannian Newton method
on the set of all matrices of fixed dimension and rank. By
``closed-form'', we mean that, besides calling an oracle for Euclidean
first and second derivatives, the method only needs to perform elementary
matrix operations, solve linear systems of equations, and compute (small-size)
matrix exponentials. By ``purely Riemannian'', we mean that it uses
the tools provided by Riemannian geometry, namely, the Riemannian
connection (instead of any other affine connection) and the Riemannian
exponential (instead of any other retraction). 

The developments strongly rely on the theory of Riemannian
submersions and are based on factorizations of low rank matrices $X$ as
$MN^\T$, where one of the factors is orthonormal. Relaxing the
orthonormality constraint is more appealing for its symmetry (the two
factors are treated alike), but it did not allow us to obtain a
closed-form expression for the Riemannian exponential.

\section*{Acknowledgement}

The authors are grateful to Nicolas Boumal for several helpful comments on a preliminary version of this work.

{\footnotesize
\bibliographystyle{alpha}
\bibliography{pabib}
}

\end{document}